\documentclass[11pt]{article}

\usepackage{amscd,amsmath, amssymb, fancyhdr, url}

\usepackage{fourier, libertine}

\newcommand{\version}{version 1.4,\ \ February 17, 2017}
\newcommand{\f}{\varphi}

\numberwithin{equation}{section}

\def\eqref#1{(\ref{#1})}
\newcommand{\goth}{\mathfrak}
\newcommand{\g}{{\mathfrak g}}

\newcommand{\Z}{{\mathbb Z}}
\newcommand{\C}{{\mathbb C}}
\newcommand{\R}{{\mathbb R}}

\def\1{\sqrt{-1}\:}
\newcommand{\restrict}[1]{{\left|_{{\phantom{|}\!\!}_{#1}}\right.}}
\newcommand{\cntrct}                % contraction with a vector field
{\hspace{2pt}\raisebox{1pt}{\text{$\lrcorner$}}\hspace{2pt}}
\newcommand{\arrow}{{\:\longrightarrow\:}}

\newcommand{\calo}{{\cal O}}
\newcommand{\cac}{{\cal C}}

\renewcommand{\bar}{\overline}
\renewcommand{\phi}{\varphi}
\renewcommand{\epsilon}{\varepsilon}
\renewcommand{\geq}{\geqslant}
\renewcommand{\leq}{\leqslant}

\newcommand{\im}{\operatorname{im}}

\newcommand{\Tot}{\operatorname{Tot}}
\newcommand{\Id}{\operatorname{Id}}

\newcommand{\const}{\operatorname{\text{\sf const}}}

\newcommand{\Spec}{\operatorname{Spec}}
\newcommand{\GL}{\operatorname{GL}}

\renewcommand{\Im}{\operatorname{Im}}

\newcounter{Mycounter}[section]
\newcounter{lemma}[section]
\newcounter{claim}[section]
\newcounter{sublemma}[section]
\newcounter{corollary}[section]
\newcounter{theorem}[section]
\renewcommand{\thetheorem}{{Theorem \thesection.\arabic{theorem}}}
\newcommand{\theorem}{%
     \setcounter{theorem}{\value{Mycounter}}
     \refstepcounter{theorem}
     \stepcounter{Mycounter}
     {\noindent \bf \thetheorem:\ }}

\newcounter{conjecture}[section]
\renewcommand{\theconjecture}{{Conjecture \thesection.\arabic{conjecture}}}
\newcommand{\conjecture}{%
     \setcounter{conjecture}{\value{Mycounter}}
     \refstepcounter{conjecture}
     \stepcounter{Mycounter}
     {\noindent \bf \theconjecture:\ }}

\newcounter{proposition}[section]
\newcounter{definition}[section]
\renewcommand{\thedefinition} {{Definition~\thesection.\arabic{definition}}}
\newcommand{\definition}{%
     \setcounter{definition}{\value{Mycounter}}
     \refstepcounter{definition}
     \stepcounter{Mycounter}
     {\noindent \bf \thedefinition:\ }}

\newcounter{example}[section]
\newcounter{remark}[section]
\renewcommand{\theremark}{{Remark \thesection.\arabic{remark}}}
\newcommand{\remark}{%
     \setcounter{remark}{\value{Mycounter}}
     \refstepcounter{remark}
     \stepcounter{Mycounter}
     {\noindent \bf \theremark:\ }}

\newcounter{problem}[section]
\newcounter{question}[section]
\makeatletter

\@addtoreset{equation}{section}
\@addtoreset{footnote}{section}
\makeatother

\def\blacksquare{\hbox{\vrule width 5pt height 5pt depth 0pt}}
\def\endproof{\blacksquare}

\newcommand{\proof}{\smallskip\noindent{\bf Proof: \ }}

\begin{document}

\begin{center}
{\LARGE\bf
Embedding of LCK manifolds \\[3mm]  with potential 
 into Hopf manifolds \\[3mm] using Riesz-Schauder theorem
}\\[5mm]
%%%%%%%%%%%%%%%%%%%%%%%%%%%%%%%%%%%%%%%%%%%%%%%%%%%%%%%%%%%%
{\large
Liviu Ornea,\footnote{Partially supported by University of Bucharest grant 1/2012.}
 and 
Misha
Verbitsky\footnote{Partially supported 
by RSCF grant 14-21-00053 within AG Laboratory NRU-HSE.\\[1mm]

}\\[4mm]
}

\end{center}

{\small
\hspace{0.15\linewidth}
\begin{minipage}[t]{0.7\linewidth}
{\bf Abstract} \\ 
A locally conformally K\"ahler (LCK) manifold with potential
is a complex manifold with a cover which admits a positive 
automorphic K\"ahler potential.  A compact LCK manifold with potential
can be embedded into a Hopf manifold,
if its dimension is at least 3. We give a functional-analytic
proof of this result based on Riesz-Schauder theorem and
Montel theorem. We provide an alternative argument for
compact complex surfaces, deducing the embedding theorem from the Spherical
Shell Conjecture.

\end{minipage}
}
%%%%%%%%%%%%%%%%%%%%%%%%%%%%%%%%%%%%%%%%%%%%%%%%

\tableofcontents

\section{Introduction}

A {\bf locally conformally K\"ahler (LCK)} manifold is a Hermitian manifold $(M,J,g)$ such that the  fundamental two-form $\omega(X,Y)=g(X,JY)$ satisfies the equation 
$$d\omega=\theta\wedge\omega$$
for a {\bf closed} one-form $\theta$, see \cite{do}. 

The one-form $\theta$ is called the {\bf Lee form}, and it produces a twisted cohomology associated to the operator $d_\theta:=d-\theta\wedge$.

An equivalent definition requires the existence of a covering $\Gamma\rightarrow \tilde M\rightarrow M$ endowed with a K\"ahler metric $\tilde g$ with respect to which the deck group $\Gamma$ acts by holomorphic homotheties. This gives rise to a character $\chi:\Gamma \rightarrow \R^{>0}$ which associates to a homothety $\gamma\in \Gamma$ its scale factor $c_\gamma$. A differential form $\eta$ on $\tilde M$ is then called {\bf automorphic} if $\gamma^*\eta=\chi(\g)\eta$. Clearly, $\tilde \omega$ is automorphic.

An LCK manifold is called {\bf with potential} if there
exists a  K\"ahler covering with K\"ahler form  form
having positive, automorphic  {\em global} potential
$\psi$: $\tilde\omega=dd^c\psi$. This is equivalent to the
existence of a function $\f$ on $M$ such that
$\omega=d_\theta d^c_\theta(\f)$, see
\cite{ov_jgp_16}. Note that $\f$ {\em is not} a potential
on $M$.

Important examples are the Hopf surfaces and, more generally, the {\bf linear Hopf manifolds} $(C^N\setminus\{0\})/\langle A\rangle$, where $\langle A\rangle$ is the cyclic group generated by a linear operator $A\in \GL(n,\C)$ with all eigenvalues $|\alpha_i| <1$, \cite{ov_pams_16}.

On the other hand, there exist compact complex manifolds
which admit LCK metrics, but no LCK metric with potential:
such are blow-ups of LCK manifolds (\cite{_vuli:blowups_})
and the LCK Inoue surfaces, \cite{O} (see also
\cite{ad}).

We proved in \cite{ov_jgp_09} (see also \cite{ov_jgp_16}) that if $(M, \omega,\theta)$ is a compact LCK manifold with potential, there exists another LCK structure $(\omega',\theta')$, close to $(\omega, \theta)$ in the $\cac^\infty$-topology, such that the corresponding $\tilde\omega'$ has a {\bf proper} potential, this being equivalent with the monodromy of the covering, $\Im(\chi)$, being isomorphic with $\Z$.

{\bf Vaisman manifolds} are LCK manifolds whose K\"ahler coverings are Riemannian cones over Sasakian manifolds (see \cite{do}, and \cite{bel} for the classification of Vaisman compact surfaces). All Vaisman manifolds are LCK with potential (represented by the squared K\"ahler norm of the pull-back of the Lee form). The converse is not true, as the example of non-diagonal Hopf manifolds shows, see \cite{ov_ma_10}. Still, the covering of an LCK manifold with potential is very close to being a cone:

\hfill

\begin{theorem}{(\cite{ov_ma_10})}\label{potcon}
Let $M$ be an LCK manifold with proper potential, $\dim_\C M\geq 3$, 
and $\widetilde M$ its K\"ahler $\Z$-covering.
Then the metric completion $\widetilde M_c$
{admits a structure of a complex variety,} actually Stein, 
compatible with the complex structure on
$\widetilde M \subset\widetilde M_c$. Moreover, 
$\tilde M_c\setminus \tilde M$ is just one point.
\end{theorem}	

\hfill

\remark
The same result seems to be true for $\dim_\C M=2$.
In Section \ref{_Kato_conj_Section_}, we deduce it from classification of
surfaces and the spherical shell conjecture on surfaces of
Kodaira class VII (\ref{_sphe_shell_Conjecture_}).

\hfill

The main property of an LCK manifold with 
potential is the folowing Kodaira type embedding result:

\hfill

\begin{theorem}  (\cite{ov_ma_10}) \label{main}
A compact LCK manifold with proper potential, of complex dimension at least $3$, can be holomorphically embedded in a linear Hopf manifold.
\end{theorem}

\hfill

\remark\label{dim_comm} The hypotesis $\dim_\C M>2$ in \ref{potcon} is essential in order to apply a result in \cite{andreotti_siu, rossi} (see \ref{asr} below) from which we deduce that the completion $\tilde M_c$ is Stein. Once we know that the completion is Stein, \ref{main} follows without further assumptions on the dimension.

\hfill

The aim of this note is to give new proofs of the above two theorems,  based on  applications of Montel and Riesz-Schauder theorems. This will require several notions of functional analysis  (see, {\em e.g.} \cite{kot}) that we recall for the reader's convenience. In the last section we comment on the possible validity of the result for complex surfaces.

\section{Preliminaries of functional analysis}

\subsection{Normal families of functions} 

\begin{definition}
Let $M$ be a complex manifold, 
and ${\cal F}$ a family of holomorphic functions
$f_i \in H^0(\calo_M)$.  ${\cal F}$ is called  
{\bf a normal family} if for each compact
$K\subset M$ there exists $C_K>0$ such that
for each $f\in {\cal F}$, $\sup_K |f| \leq C_K$.
\end{definition}

\hfill

\begin{lemma}
Let $M$ be a complex Hermitian manifold,
${\cal F}\subset H^0(\calo_M)$ a normal family,
and $K\subset M$ a compact subset. {Then there exists
	a number $A_K>0$ such that $\sup_K |f'| \leq A_K$.}
\end{lemma}

\proof
By contradiction, suppose there exists $x\in K$,
$v\in T_xM$, and a sequence $f_i \in {\cal F}$
such that $\lim_i |D_v f_i|=\infty$. We choose  a disk
$\Delta\stackrel j \hookrightarrow M$ with compact
closure in $M$, tangent to $v$ in $x$, 
such that $j(0)=x$. Let $w=tv$ have norm 1. 
Then $\sup_\Delta |f_i| < C_\Delta$ by the normal family
condition. {By Schwarz lemma (see \cite{jost}), this implies $|D_w f_i|<C_\Delta$.}
However, $t^{-1}\lim_i |D_w f_i| =\lim_i |D_v f_i|=\infty$, yielding a 
contradiction. \endproof

\subsection{Topologies on spaces of functions}
	
\begin{definition}
Let $C(M)$ be the space of functions on a topological space.
	The {\bf topology of uniform convergence on compacts} (also known
	as {\bf compact-open topology}, usually denoted as $\cac^0$) 
	is the topology on $C(M)$ whose  base of open sets
	is given by 
	\[ U(X, C):= \{f\in C(M)\ \ |\ \ \sup_K |f|< C\},
	\] for all
	compacts $K\subset M$ and $C>0$. 
	
	A sequence $\{f_i\}$ of functions converges
	to $f$ if it converges to $f$ uniformly on all compacts.
\end{definition}

\hfill
	
%\begin{remark}\label{convc0}
	%When $M$ is locally compact, 
	%any sequence of continuous functions converging in
	%$C^0$ {converges to a continuous function} 
%\end{remark}
	
\begin{remark} In a similar way one defines 
	{the  $\cac^0$-topology on the space of sections of a bundle.}
\end{remark}	

\hfill
	
%\subsubsection{The $\cac^1$ - topology on spaces of sections}

\begin{definition}
Let $B$ be a vector bundle on a smooth manifold
$M$, and $\nabla:\; B \arrow B\otimes \Lambda^1 M$ a connection. 
Define the {$\cac^1$-topology} on the space of sections of 
$B$ (denoted, as usual, by the same letter $B$)
as one where a sub-base of open sets is given by 
$\cac^0$-open sets on $B$ and $\nabla^{-1}(W)$,
where $W$ is an open set in $\cac^0$-topology in
$B\otimes \Lambda^1 M$.
\end{definition}

\hfill

\begin{remark}
A sequence $\{f_i\}$ converges in the $\cac^1$-topology if it converges 
uniformly on all compacts, and the first derivatives $\{f_i'\}$ also converge
uniformy on all compacts. {This can be seen as an equivalent definition
	of the  $\cac^1$-topology.}
\end{remark}

\subsection{Montel theorem for normal families}

\begin{theorem} {\bf (Montel).}\label{montel} 
Let $M$ be a complex manifold and ${\cal F}\subset H^0(\calo_M)$ 
a normal family of functions. Denote by $\bar{\cal F}$
its closure in the $\cac^0$-topology. {Then $\bar{\cal F}$
	is compact and contained in $H^0(\calo_M)$.}
\end{theorem}

\proof 
Let $\{f_i\}$ be a sequence of functions in ${\cal F}$.
By Tychonoff's theorem, for each compact $K$,
there exists a subsequence of $\{f_i\}$
which converges pointwise on a dense countable
subset $Z\subset K$. Taking a diagonal
subsequence, we find  a subsequence $\{f_{p_i}\}\subset \{f_i\}$
which converges pointwise on a dense 
countable subset $Z\subset M$. Since 
$|f'_i|$ is uniformly bounded on compacts,
the limit $f:= \lim_i f_i$ is  Lipschitz 
on all compact subsets of $M$. {It is thus  
	continuous,} because a pointwise limit of Lipschitz 
functions is again Lipschitz.

Then, since $|f'_i|$ is uniformly bounded on compacts, 
we can assume that $f'_i$ also converges pointwise in $Z$,
and $f:= \lim_i f_i$ is differentiable. Since a limit
of complex-linear operators is complex linear, 
$Df$ is complex linear, and $f$ is holomorphic.
{This implies that $\bar{\cal F}\cap H^0(\calo_M)$ 
	is compact. } \endproof

\subsection{The Banach space of holomorphic functions}

We begin by proving:

\hfill

\begin{theorem}
Let $M$ be a complex manifold, and $H^0_b(\calo_M)$ 
the space of all bounded holomorphic functions, equipped
with  the sup-norm $|f|_{\sup} := \sup_M |f|$.
 Then $H^0_b(\calo_M)$ is a Banach space.
\end{theorem}

\proof 
Let $\{f_i\}\in H^0_b(\calo_M)$ be a Cauchy sequence in the 
$\sup$-norm. {Then $\{f_i\}$ converges to a continuous
	function $f$} in the $\sup$-topology.

Since $\{f_i\}$ is a normal family, it
has a subsequence which converges in $\cac^0$-topology
to $\tilde f\in H^0(\calo_M)$, by Montel's Theorem 
(\ref{montel}). However, the {$\cac^0$-topology
	is weaker than the $\sup$-topology, hence $\tilde f=f$.}
Therefore, $f$ is holomorphic. \endproof

\subsubsection{Compact operators}

Recall that 
a  subset of a topological space is called {\bf 
	precompact} if its closure is compact. 

\hfill

\begin{definition}
Let $V, W$ be topological vector spaces, and let $\phi:\; V
\arrow W$ be a continuous linear operator. It
is called {\bf  compact} if the image of
any bounded set  is precompact.
\end{definition}

\hfill

\remark Note that the notion of {\bf bounded set} makes
sense in all topological vector spaces $V$. Indeed, a set
$K\subset V$ is called {\bf bounded} if for any 
open set $U\ni 0$, there exists a number $\lambda_U\in
\R^{>0}$ such that $\lambda_U K \subset U$.

\hfill

\begin{claim}\label{idcomp}
Let $V=H^0(\calo_M)$ be a space of holomorphic functions on a complex
manifold $M$ with $\cac^0$-topology. Then any 
	bounded subset of $V$ is precompact. In this case, the
identity map is a compact operator.
\end{claim}

\hfill

\proof This is a restatement of Montel's theorem
(\ref{montel}). \endproof

\hfill

\begin{remark} By Riesz theorem, {a closed ball in a normed
	vector space $V$ is never compact,} unless $V$ is
finite-dimensional. This means that $(H^0(\calo_M), \cac^0)$
does not admit a norm. A topological vector space where
any bounded subset is precompact is called {\bf Montel
  space.} 
\end{remark}

\subsubsection{Holomorphic contractions}

\begin{definition}
A {\bf contraction} of a manifold $M$ to a point
$x\in M$ is a continuous map $\phi:\; M \arrow M$ such that for any 
compact subset $K\subset M$ and any open set $U\ni x$, 
there exists $N>0$ such that for all
$n>N$, the map $\phi^n$ maps $K$ to $U$.
\end{definition}

\hfill

\begin{theorem}\label{less1}
Let $X$ be a complex variety, and let 
$\gamma:\; X \arrow X$ be a holomorphic contraction
such that $\gamma(X)$ is precompact. 
Consider the Banach space $V=H^0_b(\calo_X)$
with the sup-metric. {Then $\gamma^*:\; V \arrow V$
	is compact,} and {its operator norm
	$\|\gamma^*\|:= \sup_{|v|\leq 1}|\gamma^*(v)|$
	is strictly less than 1. }
\end{theorem}

\proof Let $B_C:=\{v\in V\  |
\   |v|_{\sup} \leq C\}$. Then 
\[|\gamma^* f|_{\sup}= \sup_{x\in \overline{\gamma(X)}}
|f(x)|.
\]
Therefore, {for any sequence $\{f_i\}$ converging in the 
	$\cac^0$-topology, the sequence $\{\gamma^* f_i\}$ converges
	in the $\sup$-topology.} However, $B_C$ is precompact in the 
$\cac^0$-topology, because it is a normal family.
Then $\gamma^* B_C$ is precompact in the 
$\sup$-topology.

Since $\sup_X |\gamma^* f|= \sup_{\gamma(X)} |f| \leq \sup_X |f|$,
one has $\|\gamma^*\|\leq 1$. If this inequality is not
strict, for some sequence $f_i\in B_1$
one has $\lim_i \sup_{x\in \gamma(X)} |f_i(x)|=1$.
Since $B_1$ is a normal family, $f_i$ has a subsequence
converging in $\cac^0$-topology to $f$. Then $\gamma(f_i)$
converges to $\gamma(f)$ in $\sup$-topology, giving 
$$\lim_i \sup_{x\in \gamma(X)} |f_i(x)|= \sup_{x\in  \gamma(X)} |f(x)|=1.$$
{Since, by the maximum principle, a holomorphic functions has no strict maxima, this means
	that $|f(x)| >1$ somewhere on $X$.} Then $f$ cannot
be the $\cac^0$-limit of $f_i\in B_1$. \endproof

\subsubsection{The Riesz-Schauder theorem}

The following result  is a Banach analogue of 
the usual spectral theorem for compact operators on Hilbert spaces. It will be the central piece in our argument.

\hfill

\begin{theorem}{\bf (Riesz-Schauder, \cite{conway})}
Let $F:\; V \arrow V$ be a compact operator on a Banach space.
Then for each non-zero $\mu \in \C$, there exists a sufficiently
large number $N\in \Z$ such that for each $n>N$ {one has 
	$$V= \ker(F-\mu\Id)^n \oplus \overline{\im (F-\mu\Id)^n},$$
	where $\overline{\im (F-\mu\Id)^n}$ is the closure
        of the image
of $(F-\mu\Id)^n$.}
Moreover, $\ker(F-\mu\Id)^n$ is finite-dimensional and independent on $n$.
\end{theorem}

\section{Proof of \ref{potcon}}

Recall that $\widetilde M_c$ denotes the metric completion of the $\Z$-cover $\tilde M$ of $M$.

\hfill

\begin{claim}\label{claim1}
{The complement $\widetilde M_c\backslash \widetilde M$ is just one
	point,} called {\bf the origin}. 
\end{claim}

\proof Indeed,
let $z_i=\gamma^{n_i}(x_i)$ be a sequence of points in
$\widetilde M$, with each $x_i$ in the fundamental
domain $\phi^{-1}([1,\lambda])$ of the $\Gamma=\Z$-action.
Clearly, the distance between two fundamental domains
$M_n:=\gamma^{n}\phi^{-1}([1,\lambda])=\phi^{-1}([\lambda^n,\lambda^{n+1}])$ and 
$M_{n+k+2}=\gamma^{n+k+2}\phi^{-1}([1,\lambda]) $
is written as
\begin{equation}\label{star}
 d(M_n,M_{n+k+2})= \sum_{i=0}^k \lambda^{n+i} v, 
\end{equation}
 where
$v$ is the distance between $M_0$ and $M_2$. Then,
$z_i$ may converge only if $\lim_i n_i=-\infty$
or if all $n_i$, except finitely many,
belong to the set $(p, p+1)$ for some $p$. The second case
is irrelevant, because each $M_i$ is compact,
and in the first case, $\{z_i\}$ is always a 
Cauchy sequence, as follows from \eqref{star}.
All such $\{z_i\}$  are therefore equivalent,
hence {converge to the same point in the metric
	completion.}
\endproof

\hfill

Recall now the following result in complex analysis:

\hfill

\begin{theorem} {(\cite{andreotti_siu}, \cite{rossi})}\label{asr} 
Let $S$ be a compact strictly pseudoconvex CR manifold, 
$\dim_\R S\geq 5$, and let $H^0(\calo_S)_b$ the ring of bounded
CR holomorphic functions. {Then $S$ is the boundary of a
Stein manifold $M$ with isolated singularities,} such that
$H^0(\calo_S)_b = H^0(\calo_M)_b$, where $H^0(\calo_M)_b$
denotes the ring of bounded holomorphic functions.
Moreover, $M$ is defined uniquely, $M=\Spec(H^0(\calo_S)_b)$.
\end{theorem}

\hfill

The proof of \ref{potcon} now goes as follows:

\hfill

\noindent{\bf  Step 1:} 
Applying Rossi-Andreotti-Siu \ref{asr}  to 
$\phi^{-1}([a, \infty[)$, we obtain a Stein variety
$\widetilde M_a$ containing $\phi^{-1}([a, \infty[)$.
Since $\widetilde M_a$ contains $\phi^{-1}([a_1, \infty[)$
for any $a_1>a$, and the  Rossi-Andreotti-Siu variety
is unique, one has $\widetilde M_a=\widetilde M_{a_1}$.
This implies that {$\widetilde M_a=: \widetilde M_c$
	is independent from the choice of $a\in \R^{>0}$. }

{It remains to identify $\widetilde M_c$ 
	with the metric completion of $\widetilde M$.}
By \ref{claim1}, {this is equivalent to the complement 
	$\widetilde M_c\backslash \widetilde M$ being a singleton.}

\smallskip

\noindent{\bf  Step 2:} The monodromy group $\Gamma=\Z$ acts on $\widetilde M_c$ by holomorphic
	automorphisms. Indeed, any holomorphic function (hence, any
holomorphic map) can be extended from $\widetilde M$ to $\widetilde M_c$
uniquely. 

\smallskip

\noindent{\bf  Step 3:} Denote by $\gamma$ the generator of $\Gamma$ which 
decreases the metric by $\lambda<1$, and let $\widetilde M^a_c$ 
be the Stein variety associated with 
$\phi^{-1}(]0, a])\subset\widetilde M$ as above. Since
$\gamma(\widetilde M^a_c)=\widetilde M^{\lambda a}_c$, for any
holomorphic function $f$ on  $\widetilde M_c$,
one has 
\[ 
\sup_{z\in\widetilde M^a_c} |f(\gamma^n(z))|= 
\sup_{z\in\widetilde M^{\lambda^n a}_c} |f(z)|\leq
\sup_{z\in\widetilde M^{a}_c} |f(z)| .
\]
{Therefore, $\{f(\gamma^n(z))\}$ is a normal family.}

\smallskip

\noindent {\bf  Step 4:}
Let $f_{\sf lim}$ be any limit point of the sequence $\{f(\gamma^n(z))\}$.
Since the sequence 
$t_i:=\sup_{z\in\widetilde M^{\lambda^i a}_c} |f(z)|$
is non-increasing, it converges, and $\sup_{z\in\widetilde M^a_c} f_{\sf lim} = \lim t_i$.
Similarly, $\sup_{z\in\widetilde M^{\lambda a}_c} f_{\sf lim} = \lim t_i$.
By the strong maximum principle, \cite{gt}, {a non-constant holomorphic function 
	on a complex manifold with boundary cannot
	have local maxima (even non-strict) outside of the boundary.}
Since $\widetilde M^{\lambda a}_c$ does not intersect the boundary of
$\widetilde M^a_c$, the function $f_{\sf lim}$ must be constant.

\smallskip

\noindent{\bf  Step 5:} Consider now the complement
$V:=\widetilde M_c\backslash \widetilde M$, and suppose it has two distinct points
$x$ and $y$. Let $f$ be a holomorphic function
which satisfy $f(x)\neq f(y)$. Replacing $f$ by an exponent of
$\mu f$ if necessarily, we may assume that $|f(x)|< |f(y)|$.
Since $\gamma$ fixes $Z$, which is compact, {for any limit $f_{\sf lim}$ 
	of the sequence $\{f(\gamma^n(z))\}$, supremum
	$f_{\sf lim}$ on $Z$ is not equal to infimum of $f_{\sf lim}$ on $Z$.}
This is impossible, hence $f=\const$ on $V$, and $V$ is one point.

This finishes the proof of \ref{potcon}.
\endproof

\section{Proof of the embedding theorem}

\newcommand{\fin}{{\text{\sf fin}}}

\ref{main} is implied by \ref{doi}. To see this, we need the following:

\hfill

\begin{definition}\label{A-finite}
Let $\gamma$ be an endomorphism of a vector space $V$.
A vector $v\in V$ is called {\bf $\gamma$-finite} if the
subspace $\langle v, \gamma(v), \gamma^2(v),\ldots \rangle$
is finite-dimensional.
\end{definition}

\hfill

\begin{theorem}\label{doi}
	 Let $M$ be an LCK manifold with
potential, $\dim_\C M>2$, and $\widetilde M$ its K\"ahler $\Z$-covering.
Consider the metric completion $\widetilde M_c$ with its complex
structure and a contraction $\gamma:\; \widetilde M_c \arrow \widetilde M_c$ 
generating the
$\Z$-action. Let $H^0(\calo_{\widetilde M_c})_\fin$ be the space of
functions which are $\gamma^*$-finite.  Then
	$H^0(\calo_{\widetilde M_c})_\fin$ is dense in the $\sup$-topology
	on each compact subset of $\widetilde M_c$.
\end{theorem}

\subsection{\ref{doi} implies \ref{main}}

{\noindent\bf Step 1:} 
Let $W\subset H^0(\calo_{\widetilde M_c})_\fin$ be an
$m$-dimensional $\gamma^*$-invariant subspace $W$ with basis
$\{w_1, \ldots, w_m\}$. Then the following diagram is
commutative:
\begin{equation*}
\begin{CD}
\widetilde M@>\Psi>> \C^m \\
@V{\gamma}VV  @VV{\gamma^*}V \\
M@>\Psi >>  \C^m
\end{CD},
\end{equation*}
where $\Psi(x)=(w_1(x), w_2(x), \ldots, w_m(x))$.

Suppose that the map $\Psi$ associated with a given 
$W\subset  H^0(\calo_{\widetilde M_c})_\fin$ is injective. Then {the
	quotient map gives an embedding 
	$\Psi:\; \widetilde M/\Z \arrow (\C^m\backslash 0)/\gamma^*$;}
all eigenvalues of $\gamma^*$ are $<1$ because its operator
norm is $<1$, by \ref{less1}. 

\smallskip

{\noindent\bf Step 2:} To find an appropriate 
$W\subset  H^0(\calo_{\widetilde M_c})_\fin$, choose a holomorphic embedding
$\Psi_1:\; \widetilde M_c\hookrightarrow \C^n$, which exists because
$\widetilde M_c$ is Stein. Let $\tilde w_1, \ldots, \tilde w_n$
be the coordinate functions of $\Psi_1$. 
\ref{doi} allows one to approximate $\tilde w_i$ by 
$w_i \in H^0(\calo_{\widetilde M_c})_\fin$ in $\cac^0$-topology.
Choosing $w_i$ sufficiently close to $\tilde w_i$
	in a compact fundamental domain of the $\Z$-action, we obtain
	that $x\mapsto (w_1(x), w_2(x), ..., w_n(x))$ is injective
	in a compact fundamental domain of $\Z$.

Finally,
take $W\subset H^0(\calo_{\widetilde M_c})_\fin$ generated by 
the $\gamma^*$ from $w_1, \ldots, w_n$, and apply Step 1.
\endproof

\smallskip

\def\m{{{\goth m}}}

\subsection{Proof of \ref{doi}} 

The core of our argument is an application of Riesz-Schauder theorem. First we prove:

\hfill

\begin{proposition}\label{propunu}
	Fix a precompact subset $\widetilde M^a_c:= \phi^{-1}([0, a[)$,
	where $\phi:\; \widetilde M_c:\arrow \R^{>0}$ is the K\"ahler potential.
	Let $A$ be the ring of bounded holomorphic functions
	on $\widetilde M^a_c$, and ${\goth m}$ the maximal ideal of the
	origin point. Clearly, $\gamma^*$ preserves $\m$ and all 
	its powers. Let $P_k(t)$ be the minimal polynomial
	of $\gamma^*\restrict{A/\m^k}$. Then $\im(P_k(\gamma^*))\subset \m^k(A)$,
		and $\ker P_k(\gamma^*)$ generates $A/\m^k$.
\end{proposition}

\proof Since $P_k(t)$ is a minimal polynomial of
$\gamma^*$ on $A/\m^k$, the endomorphism $P_k(\gamma^*)$
acts trivially on $A/\m^k$, by Cayley-Hamilton theorem, hence it maps $A$ to $\m^k$.

From Riesz-Schauder theorem applied to the Banach space
$A$ and $F=P_k(\gamma^*)-P_k(0)$, with $\mu=-P_k(0)$,  
it follows that $A=\ker(P_k(\gamma^*)\oplus  \overline{\im (P_k(\g^*))^n}$. Since $P_k(\gamma^*)$ acts trivially on $A/\m^k$, its image lies in $\m^k$. This gives a surjection
of $\ker P_k(\gamma^*)$ onto $A/\m^k$.
\endproof

\hfill

This implies:

\hfill

\begin{proposition}\label{density}
Let $H^0(\calo_{\widetilde M_c})_\fin\subset H^0(\calo_{\widetilde M_c})$
be the set of $\gamma^*$-finite functions and $\m$
the maximal ideal of the origin in $\widetilde M_c$. Then
	$H^0(\calo_{\widetilde M_c})_\fin$ is dense in $\m$-adic topology\footnote{Recall that for a ring $A$ with a proper  ideal  $\m$, 
		the $\m$-adic topology on $A$ is given by the base of open sets formed by  $\m^k$ and their translates.}.
\end{proposition}

\proof A subspace $V\subset A$ is dense in $\m$-adic topology
in $A$ $\Leftrightarrow$ the quotient $V/v\cap \m^k$
surjects to $A/\m^k$. This is proven in\ref{propunu} 
for the ring of bounded holomorphic functions
on $\widetilde M^a_c$. However, any such function can be
extended to $\gamma^*$-finite function on $\widetilde M_c$
using the $\gamma^*$-action. \endproof

\hfill

To finish the proof, observe that \ref{doi} is implied by the following:

\hfill

\begin{claim}
Let $X$ be a connected complex variety,  
$A$ the ring of bounded holomorphic functions on $X$, $x\in X$ a point,
$\m\subset A$ its maximal ideal, and $R:\; A \arrow \hat A$
the natural map from $A$ to its $\m$-adic completion. Then $R$ is continuous in $\cac^0$-topology and induces
	homeomorphism of any bounded set  to its image.
\end{claim}

\proof Continuity is clear because the $\cac^0$-topology on holomorphic functions 
is equivalent to $\cac^1$-topology, $\cac^2$-topology and
so on, by Montel Theorem (\ref{montel}). 
Therefore, taking successive derivatives in a point is continuous in
$\cac^0$-topology. However, $R$ takes a function and replaces
it by its Taylor series.

To see that $R$ is a homeomorphism, notice that any
bounded, closed subset of $A$ is compact, hence its image
under a continuous map is also closed. Then $R$ induces
a homeomorphism on all bounded sets. To see that the  
preimage of a converging sequence is converging, notice
that any such sequence is bounded in $A$ by another appication
of Schwarz lemma.
\endproof

%%%%%%%%%%%%%%%%%%%%%%%%%%%%%%%%%%%%%%%%%%%%%%%%%%%%%%%%%%%%%%%%

\section{Kato conjecture and non-K\"ahler surfaces}
\label{_Kato_conj_Section_}

%%%%%%%%%%%%%%%%%%%%%%%%%%%%%%%%%%%%%%%%%%%%%%%%%%%%%%%%%%%%%%%%

All surfaces in this section are assumed to be compact.

\hfill

\definition
A complex surface $M$ with $b_1(M)=1$ and Kodaira dimension $-\infty$
is called {\bf Kodaira class VII surface}.
If it is also minimal, it is called {\bf class VII${}_0$ surface}.

\hfill

\remark
Class VII surfaces are obviously non-K\"ahler. Indeed, their $b_1$ is odd.

\hfill

The main open question in the classification of non-K\"ahler surfaces is
the following conjecture, called {\bf spherical shell conjecture},
or {\bf Kato conjecture}. To state it, recall first:

\hfill

%%%%%%%%%%%%%%%%%%%%%%%%%%%%%%%%%%%%%%%%%%%%%%%%%%%%%%%%%%%%%%%%%%%%
\definition
Let $S\subset M$ be a real submanifold  in a complex
surface, diffeomorphic to $S^3$.
We call $S$ {\bf sphericall shell} if $M\backslash S$ is connected,
and $S$ has a neighbourhood which is biholomorphic to an annulus in 
$\C^2$. A class VII${}_0$ surface which contains a spherical
shell is called a {\bf Kato surface}.

\hfill

\remark
From \cite{kato}, we know that any Kato surface contains
exactly $b_2(M)$ distinct rational curves (the converse was proven in \cite{_DOT:Kato_surfaces_}).

\hfill

%%%%%%%%%%%%%%%%%%%%%%%%%%%%%%%%%%%%%%%%%%%%%%%%%%%%%%%%%%%%%
\conjecture\label{_sphe_shell_Conjecture_} 
(spherical shell conjecture) 
Any class VII${}_0$ surface with $b_2>0$ is a Kato surface.

\hfill

%%%%%%%%%%%%%%%%%%%%%%%%%%%%%%%%%%%%%%%%%%%%%%%%%%%%%%%%
\theorem\label{_sphe_implies_dim2_Theorem_}
Assume that the spherical shell conjecture is true.
Then \ref{potcon} and \ref{main} are true in dimension 2.

\hfill

\noindent{\bf Proof:} The only part of the proof missing for
dimension 2 is Rossi-Andreotti-Siu Theorem (\ref{asr}), see also \ref{dim_comm}.
We used it to prove the following result (which is
stated here as a conjecture, because we don't know
how to prove it for class VII${}_0$ non-Kato surfaces).

\hfill

%%%%%%%%%%%%%%%%%%%%%%%%%%%%%%%%%%%%%%%%%%%%%%%%%%%%
\conjecture\label{_pot_then_ras_Conjecture_}
Let $M$ be an LCK complex surface with proper
potential, and $\tilde M$ its K\"ahler $\Z$-cover.
Then the metric completion of $\tilde M$, realized by adding just one point,  is a
Stein variety.

\hfill

%%%%%%%%%%%%%%%%%%%%%%%%%%%%%%%%%%%%%%%%%%%%%%%%%%%%
\remark
For $\dim M\geq 3$, this is \ref{potcon}

\hfill

\ref{_pot_then_ras_Conjecture_}
follows from the spherical shell conjecture and the
classification of surfaces.

\hfill

First of all, notice that an LCK surface $M$ with an LCK potential
cannot contain rational curves. Indeed, if $M$ contains
rational curves, by homotopy lifting $\tilde M$ would also
contain rational curves, but $\tilde M$ is embedded
to a Stein variety. This implies that $M$ cannot
be a Kato surface, and that $M$ is minimal.

\hfill

If the spherical shell conjecture is true,
class VII surfaces which are not Kato have $b_2=0$.
However, class VII surfaces with $b_2=0$ were
classified by Bogomolov, Li, Yau, Zheng and Teleman
(\cite{_Bogomolov:VII_,_Bogomolov:VII_,_LYZ:flat_,_Li_Yau_,_Teleman:bogomolov_,_Teleman:donaldson_}).
From these works it follows that any class VII surface with $b_2=0$
is biholomorphic to a Hopf surface or to an Inoue surface.

\hfill

Inoue surfaces don't have LCK potential for topological
reasons (\cite{O}). The Hopf surfaces are quotients
of $\C^2\setminus 0$ by an action of $\Z$, hence their 1-point completions
are affine, and hence Stein.

\hfill

The only non-K\"ahler minimal surfaces which are not of class
VII are non-K\"ahler elliptic surfaces
(\cite{_Barth_Peters_Van_de_Ven_}).
These surfaces are obtained as follows.
Let $X$ be a 1-dimensional compact complex orbifold,
and $L$ an ample line bundle on $X$. Consider
the space $\tilde M$ of all non-zero vectors in the
total space of $L^*$, and let $\Z$ act on $\tilde M$
as $v \mapsto \alpha v$, where $\alpha\in \C$
is a fixed complex number, $|\alpha|>1$. 
Any non-K\"ahler elliptic surface is isomorphic
to $\tilde M/\Z$ for appropriate $\alpha$, $X$ and $L$.
However, the sections of $L^{\otimes n}$ 
define holomorphic functions on $\tilde M\subset
\Tot(L^*)$, identifying $\tilde M$ 
and the corresponding cone over $X$. As such, $M$ is Vaisman (\cite{bel}), in particular LCK with potential. The completion of this
cone is $\tilde M_c$, and it is affine
(\cite[\S 8]{_EGA2_}), and hence Stein.
This finishes the proof of \ref{_sphe_implies_dim2_Theorem_}.
\endproof

\hfill

\noindent{\bf Acknowledgment:} The authors thank Georges Dloussky for his kind advice and for a bibliographical information, and the anonymous referee for very useful remarks.

\hfill

{\small

\noindent {\sc Liviu Ornea\\
University of Bucharest, Faculty of Mathematics, \\14
Academiei str., 70109 Bucharest, Romania}, and:\\
{\sc Institute of Mathematics "Simion Stoilow" of the Romanian
Academy,\\
21, Calea Grivitei Str.
010702-Bucharest, Romania\\
\tt lornea@fmi.unibuc.ro, \ \  Liviu.Ornea@imar.ro}

\hfill

\noindent {\sc Misha Verbitsky\\
Laboratory of Algebraic Geometry, \\
Faculty of Mathematics, National Research University HSE,\\
7 Vavilova Str. Moscow, Russia}, also: \\
{\sc Universit\'e Libre de Bruxelles, D\'epartement de Math\'ematique\\
Campus de la Plaine, C.P. 218/01, Boulevard du Triomphe\\
B-1050 Brussels, Belgium\\
\tt verbit@verbit.ru, mverbits@ulb.ac.be }
}

\begin{thebibliography}{100}

\bibitem[AD]{ad} V. Apostolov, G. Dloussky, {\em On the Lee classes of locally conformally symplectic complex surfaces}, arXiv:1611.00074.

\bibitem[AS]{andreotti_siu} A. Andreotti, Y.T. Siu,  {\em Projective embeddings of pseudoconcave spaces}, Ann. Scuola Norm. Sup. Pisa {\bf 24}, 231--278 (1970). 

\bibitem[BHPV]{_Barth_Peters_Van_de_Ven_}
Barth W., Hulek K., Peters C., Van de Ven A.
{\em Compact complex surfaces}, 2004, Springer Verlag.


\bibitem[Be]{bel} F.A. Belgun, {\em On the metric structure of non-K\"ahler complex surfaces}, Math. Ann. {\bf 317} (2000), 1--40.


\bibitem[Bo1]{_Bogomolov:VII_}  
F. A. Bogomolov, {\em Classification of surfaces of
class VII${}_0$ with $b_2(M)=0$,} Math. USSR Izv., vol. 10,
pp. 255-269, 1976.

\bibitem[Bo2]{_Bogomolov:affine_} 
F. A. Bogomolov, {\em Surfaces of class     and affine
geometry,} Izv. Akad. Nauk SSSR Ser. Mat., vol. 46,
iss. 4, pp. 710-761, 896, 1982.

\bibitem[Co]{conway} J. B. Conway, A course in functional analysis, G.T.M. {\bf 96}, Springer 1990.

\bibitem[DOT]{_DOT:Kato_surfaces_}
Dloussky, Georges; Oeljeklaus, Karl; Toma, Matei (2003),
{\em Class VII0 surfaces with b2 curves}, The Tohoku
Mathematical Journal. Second Series, 55 (2): 283-309,

\bibitem[EGA2]{_EGA2_}
Grothendieck, Alexandre; Dieudonn\'e, Jean
 {\em \`El\'ements de g\'eom\'etrie alg\'ebrique: II. \'Etude
 globale \'el\'ementaire de quelques classes de
 morphismes.} Publications Math\'ematiques de l'IH\'ES. 8 (1961). 



\bibitem[DO]{do} S. Dragomir, L. Ornea, Locally conformally K\"ahler manifolds, Progress in Math. {\bf 55}, Birkh\"auser, 1998.

\bibitem[GT]{gt} D. Gilbarg, N.S. Trudinger, Elliptic partial differential equations of second order,  Classics in Mathematics. Springer-Verlag, Berlin, 2001.

\bibitem[Jo]{jost} J. Jost, Compact Riemann surfaces, Springer, 2002.

\bibitem[Kat]{kato} Ma. Kato, {\em Compact complex manifolds containing \\global'' spherical shells}, Proc. Japan Acad. {\bf 53} (1977), no. 1, 15--16.

\bibitem[Kot]{kot} G. K\"othe, Topological vector spaces, I, II, Springer-Verlag, 1969, 1979. 

\bibitem[LY]{_Li_Yau_} 
J. Li and S. Yau, 
{\em Hermitian-Yang-Mills connection on non-K\"ahler
  manifolds,} in Mathematical Aspects of String Theory,
Singapore: World Sci. Publishing, 1987, vol. 1,
pp. 560-573.

\bibitem[LYZ]{_LYZ:flat_}
J. Li, S. Yau, and F. Zheng, {\em On projectively flat
Hermitian manifolds,} Comm. Anal. Geom., vol. 2, iss. 1,
pp. 103-109, 1994.


\bibitem[OV1]{ov_jgp_09} L. Ornea, M. Verbitsky, {\em Morse-Novikov cohomology of locally conformally K\"ahler manifolds}, J. Geom. Phys. {\bf 59}, No. 3 (2009), 295--305.

\bibitem[OV2]{ov_ma_10} L. Ornea, M. Verbitsky, {\em Locally conformal K\"ahler manifolds with potential}, Mathematishe Annalen {\bf 348} (2010), 25--33.

\bibitem[OV3]{ov_jgp_16} L. Ornea, M. Verbitsky, {\em LCK rank of locally conformally K\"ahler manifolds with potential}, J. Geom. Phys. {\bf 107} (2016), 92--98.

\bibitem[OV4]{ov_pams_16} L. Ornea, M. Verbitsky, {\em Locally conformally K\"ahler metrics obtained from pseudoconvex shells}, Proc. Amer. Math. Soc. {\bf 144} (2016), 325--335. 

\bibitem[Ot]{O} 
A. Otiman, {\em Morse-Novikov cohomology of locally
  conformally K\"ahler surfaces},  arXiv:1609.07675.


\bibitem[Ro]{rossi} H. Rossi, {\em Attaching analytic spaces to an analytic space along a pseudo--convex boundary}, Proceedings of the Conference Complex Manifolds (Minneapolis), pp. 242--256. Springer, Berlin (1965).


\bibitem[Te1]{_Teleman:bogomolov_}
A. D. Teleman, 
{\em Projectively flat surfaces and Bogomolov's theorem on
  class VII${}_0$ 
surfaces,} Internat. J. Math., vol. 5, iss. 2, pp. 253-264, 1994.

\bibitem[Te2]{_Teleman:donaldson_}
A. D. Teleman, {\em Donaldson theory on
non-K\"ahlerian surfaces and class VII surfaces with $b_2=0$,}
Invent. Math., vol. 162, iss. 3, pp. 493-521, 2005.



\bibitem[Vu]{_vuli:blowups_} 
V. Vuletescu, {\em Blowing-up points on locally
conformally
K\"ahler manifolds},
Bull. Math. Soc. Sci. Math. Roumanie {\bf 52}(100) (2009), 387--390.


\end{thebibliography}
\end{document}